\font\gothic=eufm10.
\font\sets=msbm10.
\font\stampatello=cmcsc10.
.

\def\1{{\bf 1}}
\def\square{\hbox{\vrule\vbox{\hrule\phantom{s}\hrule}\vrule}}
\def\defineq{\buildrel{def}\over{=}}

\def\C{\hbox{\sets C}}
\def\N{\hbox{\sets N}}
\def\P{\hbox{\sets P}}

\def\Z{\hbox{\sets Z}}
\def\SingSer{\hbox{\gothic S}}

\par
\centerline{\bf An elementary property of correlations}
\bigskip
\centerline{\stampatello giovanni coppola}\footnote{ }{MSC $2010$: $11{\rm N}05$, $11{\rm P}32$, $11{\rm N}37$ - Keywords: correlation, shift Ramanujan expansion, $2k-$twin primes} 
\bigskip
\bigskip
\bigskip
\par
\centerline{\stampatello 1. Introduction and statement of the results.}
\smallskip
\par
\noindent
We define for any arithmetic functions $f,g:\N \rightarrow \C$ their {\it correlation} (or shifted convolution sum) of {\it shift} $a$: 
$$
C_{f,g}(N,a)\defineq \sum_{n\le N}f(n)g(n+a),
\enspace
\forall a\in \N. 
$$
\par
\noindent
Notice in passing that it is an arithmetic function itself, of argument $a\in \N$, the shift. In fact, in $\S5$ of [CMu2] we introduced the {\it shift-Ramanujan expansion}, i.e. (see $(1)$ in [CMu2] for \thinspace $c_{\ell}(a)$, the {\it Ramanujan sum}) : 
$$
C_{f,g}(N,a)=\sum_{\ell=1}^{\infty}\widehat{C_{f,g}}(N,\ell)c_{\ell}(a),
\enspace
\forall a\in \N. 
$$
\par
\noindent
Any arithmetic function $F:\N \rightarrow \C$ may be written as \thinspace $F(n)=\sum_{d|n}F'(d)$, by M\"{o}bius inversion [T], with a uniquely determined $F'\defineq F\ast \mu$ (see [T] for $\ast$, $\mu$), its {\it Eratosthenes transform} (Wintner's [W] terminology). 
\par
We shall, hereafter, truncate \thinspace $g(m)=\sum_{q|m}g'(q)$ \thinspace as \thinspace $g_N(m)\defineq \sum_{q|m,q\le N}g'(q)$; in fact, our calculations will be shorter, with an $a-$independent truncation at a small cost, i.e. the error is small: 
$$
C_{f,g}(N,a)-C_{f,g_N}(N,a)=\sum_{N<q\le N+a}g'(q)\sum_{{n\le N}\atop {n\equiv -a\bmod q}}f(n)
\ll \max_{n\le N}|f(n)| \cdot \max_{N<q\le N+a}|g'(q)|\cdot a, 
\enspace
\forall a\in \N, 
\leqno{(1)}
$$
\par
\noindent
which, in the case $f$ and $g$ satisfy the Ramanujan Conjecture\footnote{$^1$}{\noindent Ramanujan Conjecture for $f$ says: $f(n)\ll_{\varepsilon} n^{\varepsilon}$, as $n\to \infty$. Hereafter Vinogradov's $\ll$ is equivalent to Landau's $O-$notation, [T], also, $\ll_{\varepsilon}$ says, like $O_{\varepsilon}$, that the constant may depend on arbitrarily small $\varepsilon>0$.}, is \enspace $O_{\varepsilon}\left(N^{\varepsilon}\left(N+a\right)^{\varepsilon}a\right)$, uniformly \enspace $\forall a\in \N$. 
\par
We say, by definition, that a correlation \enspace $C_{f,g}(N,a)$ \enspace is {\it fair} when the dependence on the shift $a$ is only inside the argument of $g$, $n+a$, but not in $f$, $g$, neither in their supports. Assuming \lq \lq $g$ {\it has range} $Q$\rq \rq, i.e., 
$$
g(m)=g_Q(m)\defineq \sum_{q|m,q\le Q}g'(q)=\sum_{\ell \le Q}\hat{g}(\ell)c_{\ell}(m),
\enspace 
\hbox{\rm where}
\enspace 
\hat{g}(\ell)\defineq \sum_{q\equiv 0\bmod \ell}{{g'(q)}\over q}
$$
\par
\noindent
(that is, compare [CMu1], $g_Q$ finite Ramanujan expansion), with $Q$ independent of $a$, then \enspace $C_{f,g}(N,a)$ \enspace is  
$$
C_{f,g}(N,a)=C_{f,g_Q}(N,a)=\sum_{q\le Q}\hat{g}(q)\sum_{n\le N}f(n)c_q(n+a), 
\enspace
\forall a\in \N, 
\leqno{(2)}
$$
\par
\noindent
where the $\hat{g}(q)$ are above Ramanujan coefficients of $g$. This correlation is fair iff (i.e., \lq \lq if and only if\rq \rq) all the $f(n)$, the $\hat{g}(q)$ \& their supports don't depend on $a$, i.e.: $a-$dependence is only in $c_q(n+a)$ ! 
We define: 
$$
C'_{f,g_N}(N,\ell)\defineq \sum_{t|\ell}C_{f,g_N}(N,t)\mu\left({{\ell}\over t}\right), 
$$
\par
\noindent
which has part in the following Delange Hypothesis (DH,here), for the truncated correlation $C_{f,g_N}(N,a)$: 
$$
\sum_{d=1}^{\infty}{{2^{\omega(d)}}\over d}\left|C'_{f,g_N}(N,d)\right|<\infty, 
\leqno{\rm (DH)}
$$
\par				
\noindent
where the arithmetic function $\omega(d)$ counts the prime factors of $d$, whence $2^{\omega(d)}$ is the number of square-free divisors of $d$, that has bound 
$$
2^{\omega(d)}\ll_{\varepsilon} d^{\varepsilon},
\enspace
\hbox{\rm as}
\enspace 
d\to \infty, 
$$
\par
\noindent
since it is bounded by the number of divisors of $d$ (and divisor function also satisfies Ramanujan Conjecture). 
\par
The ones listening to my talk of 5 Sep 2017, in Poznan, Poland, at NTW2017 (see on ResearchGate) will remember, probably, that $(DH)$ implies Carmichael's Formula (in general, see the following): here 
$$
\widehat{C_{f,g_N}}(N,\ell)={1\over {\varphi(\ell)}}\lim_{x\to \infty}{1\over x}\sum_{a\le x}C_{f,g_N}(N,a)c_{\ell}(a), 
\leqno{\rm (CF)}
$$
\par
\noindent
where \enspace $\varphi(\ell)\defineq \left|\left\{n\le \ell : (n,\ell)=1\right\}\right|$ is the Euler function. Actually, the implication $(DH)\Rightarrow (CF)$ follows from a result of Wintner (of 1943 [W]) and a result of Delange (published in 1976, [De]) that we quote here from [ScSp] Theorem 2.1 in Chapter VIII on Ramanujan expansions (restating and selecting properties), for all arithmetic functions $F$ : 
\medskip
\par
\noindent {\bf Wintner-Delange Formula}. {\it Let } $F:\N \rightarrow \C$ {\it satisfy Delange Hypothesis, namely}
$$
\sum_{d=1}^{\infty}{{2^{\omega(d)}}\over d}\, \left|F'(d)\right|<\infty. 
$$
\par
\noindent
{\it Then the Ramanujan expansion}
$$
\sum_{q=1}^{\infty}\widehat{F}(q)c_q(n)
$$ 
\par
\noindent
{\it converges pointwise to } $F(n)$, $\forall n\in \N$, {\it with coefficients given by the formula}
$$
\widehat{F}(q)=\sum_{d\equiv 0\bmod q}{{F'(d)}\over d},
\enspace
\forall q\in \N
$$ 
\par
\noindent
{\it (where the series on RHS, right hand side, converges pointwise, $\forall q\in \N$) and also by Carmichael\footnote{$^2$}{\rm The name given here is in honour of Carmichael [Ca] : maybe, compare [Mu, pp.26-27], it's Wintner's} formula}
$$
\widehat{F}(q)={1\over {\varphi(q)}}\lim_{x\to \infty}{1\over x}\sum_{n\le x}F(n)c_q(n), 
\enspace
\forall q\in \N
$$ 
\par
\noindent
({\it where the limit on RHS exists in complex numbers}, $\forall q\in \N$)
\medskip
\par
\noindent
We don't need, actually, to prove this result, as it follows from (quoted) Th.2.1 of [ScSp]. In the case $F(a)=C_{f,g_N}(N,a)$, assuming the above $(DH)$ (i.e., Delange Hypothesis for present $F$), then Wintner-Delange formula implies the above $(CF)$ (i.e., Carmichael Formula for $F$); this, in turn, is condition $(ii)$ of Theorem 1 in [Cmu2] which is equivalent, choosing $Q=N$, to the following {\it {\bf R}amanujan {\bf e}xact {\bf e}xplicit {\bf f}ormula} (as I named condition $(iii)$ in Theorem 1 [CMu2]) for $C_{f,g_N}$, that is also uniform in $a\in \N$ : 
$$
C_{f,g_N}(N,a)=\sum_{\ell \le N}\left({{\hat{g}(\ell)}\over {\varphi(\ell)}}\sum_{n\le N}f(n)c_{\ell}(n)\right)c_{\ell}(a). 
\leqno{\hbox{\stampatello R.e.e.f.}}
$$
\par
\noindent
This part of our original correlation $C_{f,g}$, for general $f,g:\N \rightarrow \C$ satisfying Ramanujan Conjecture, has a lot of structure (it's a truncated divisor sum!); adding the other part, we estimated above in $(1)$, we get, for fair correlations with $(DH)$, the following \lq \lq structure $+$ small error\rq \rq$-$elementary property (that gives name to the paper). 
\smallskip
\par				
\noindent {\bf Theorem.} {\it Let } $f,g:\N \rightarrow \C$ {\it satisfy the Ramanujan Conjecture and be such that, for the $N-$truncated divisor sum } $g_N(m)$ {\it defined above, the correlation } $C_{f,g_N}$ {\it is fair and satisfies $(DH)$. Then } 
$$ 
C_{f,g}(N,a)=\sum_{\ell \le N}\left({{\hat{g}(\ell)}\over {\varphi(\ell)}}\sum_{n\le N}f(n)c_{\ell}(n)\right)c_{\ell}(a)
 +O_{\varepsilon}\left(N^{\varepsilon}\left(N+a\right)^{\varepsilon}a\right),
$$
\par
\noindent
{\it uniformly in } $a\in \N$. 
\medskip
\par
\noindent
What said up to now suffices to prove the Theorem (notice: $(1)$ \& $(2)$, Wintner-Delange result above and Theorem 1 in [CMu2] are the whole {\bf proof} ). QED\footnote{$^3$}{In this paper, QED($=$Quod Erat Demonstrandum$=$What was to be shown) is not the end of the story, in a proof (we use $\square$ for it); also, in the following, it will indicate an involved, smaller, part of proof ending} 
\par
However, thanks to the importance and generality (in $\S3$ we have, say, a huge application too) we will provide a step-by-step Proof in next section, $\S2$. 
\medskip
\par
In a perfectly similar fashion to the Proof of Corollary 1 [CMu2], from Theorem 1 [CMu2], we can prove (but we will not do) the following consequence.
\smallskip
\par
\noindent {\bf Corollary.} {\it Assume } $f,g:\N \rightarrow \C$ {\it satisfy Ramanujan Conjecture, where furthermore } $f$ {\it is a $D-$truncated divisor sum, say } $f(n)=f_D(n)\defineq {\displaystyle \sum_{d|n,d\le D} }f'(d)$, {\it with } ${{\log D}\over {\log N}}<1-\delta$. {\it Also, let the  correlation } $C_{f,g_N}$ {\it be fair, with } $(DH)$. {\it Then } 
$$
C_{f,g}(N,a)=\SingSer_{f,g}(a)N+O\left(N^{1-\delta}\right)+O_{\varepsilon}\left(N^{\varepsilon}\left(N+a\right)^{\varepsilon}a\right),
$$
\par
\noindent
{\it uniformly in } $a\in \N$, {\it where the, say, \lq \lq singular sum\rq \rq, here, is defined with $f,g$ Ramanujan coefficients as}
$$
\SingSer_{f,g}(a)\defineq \sum_{q\le N}\hat{f}(q)\hat{g}(q)c_q(a),
\enspace
\forall a\in \N.
$$
\medskip
\par
\noindent
Before an \lq \lq unnecessary\rq \rq, but beautiful, Proof of our Theorem (that, actually, will prove even the above Wintner-Delange formula, I mentioned in my talk), we apply our Theorem, in section $\S3$, to the noteworthy case of $2k-$twin primes, assuming $(DH)$ for them. Also, I realized later that, not like I told in the talk, this noteworthy case also comes from our Theorem 1 [CMu2]. In fact, truncating $g$ at $Q=N$ (in Theorem 1) and considering a kind of approximation, to original correlation, as given above in equation $(1)$, everthing works fine!

\bigskip
\bigskip

\par
\centerline{\stampatello 2. The detailed proof of our Theorem.}
\smallskip
\par
\noindent {\bf Proof.} Starting from $(1)$, we are left with the task of proving the {\bf Reef} above, i.e., 
$$
\sum_{n\le N}f(n)\sum_{q|n+a,q\le N}g'(q)=\sum_{\ell \le N}\left({{\hat{g}(\ell)}\over {\varphi(\ell)}}\sum_{n\le N}f(n)c_{\ell}(n)\right)c_{\ell}(a). 
$$
\par
\noindent
The hypotheses of our Theorem ensure that the LHS, namely $C_{f,g_N}(N,a)$, satisfies $(DH)$ above. Now, we need to infer $(DH)$ $\Rightarrow$ $(CF)$ (see the above), namely, get the Carmichael formula for our $C_{f,g_N}(N,a)$, so to have in the following, say, a way to infer the R.e.e.f. ! However, we'll supply even more, by providing a proof, for the above \lq \lq Wintner-Delange formula\rq \rq. (Hence, in the immediate following we'll import arguments from [De] \& [ScSp].) 
\par
In order to prove it, we wish to prove that the following double series, over $\ell,d$ summations, is absolutely convergent; so, we may write the equation expressing it in two ways (first summing over $\ell$, then $d$ and the vice versa) : 
$$
\sum_{d=1}^{\infty}\sum_{\ell|d}{{F'(d)}\over {d}}c_{\ell}(n)
=\sum_{\ell=1}^{\infty}\sum_{d\equiv 0\bmod \ell}{{F'(d)}\over {d}}c_{\ell}(n), 
\enspace 
\forall n\in \N, 
\leqno{(\ast)}
$$
\par				
\noindent
namely, exchange sums. In fact, ${\displaystyle {1\over d}\sum_{\ell|d}c_{\ell}(n) }=\1_{d|n}$, for $\1_{\wp}=1$ iff $\wp$ is true ($0$ otherwise), [CMu, Lemma 1] gives LHS 
$$
\sum_{d=1}^{\infty}{{F'(d)}\over {d}}\sum_{\ell|d}c_{\ell}(n)
=\sum_{d|n}F'(d)=F(n), 
$$
\par
\noindent
with on RHS the Wintner-Delange coefficients 
$$
\sum_{d\equiv 0\bmod \ell}{{F'(d)}\over {d}}, 
\enspace \forall \ell\in \N 
$$
\par
\noindent
thus supplying a proof of the first (Wintner-Delange's!) formula and also ensuring pointwise convergence of Ramanujan expansion, with these coefficients: 
$$
(\ast)
\enspace \Rightarrow \enspace
F(n)=\sum_{\ell=1}^{\infty}\left(\sum_{d\equiv 0\bmod \ell}{{F'(d)}\over {d}}\right)c_{\ell}(n),
\forall n\in \N. 
$$
\par
\noindent
Absolute convergence of double series comes from the fact that LHS with moduli, $\forall d,\ell\in \N$, are bounded by 
$$
\sum_{d=1}^{\infty}{{\left|F'(d)\right|}\over {d}}\sum_{\ell|d}\left|c_{\ell}(n)\right|
\le n\sum_{d=1}^{\infty}{{\left|F'(d)\right|}\over {d}}2^{\omega(d)}<\infty, 
\enspace
\forall n\in \N, 
$$
\par
\noindent
coming as we know from Delange Hypothesis, starting from the optimal bound, proved by Hubert Delange: 
$$
\sum_{\ell|d}\left|c_{\ell}(n)\right|\le n\cdot 2^{\omega(d)}, 
$$
\par
\noindent
for which we refer to Delange's original paper [De] (also, for comments about optimality). 
\par
Left to prove, for Wintner-Delange formula above, is the fact that above coefficients (Wintner-Delange's, which we know, now, to be the Ramanujan coefficients!) are given also by the Carmichael formula: 
$$
{1\over {\varphi(q)}}\lim_{x\to \infty}{1\over x}\sum_{n\le x}F(n)c_q(n)=\sum_{d\equiv 0\bmod q}{{F'(d)}\over d}, 
$$
\par
\noindent
our task, now; for which we plug (in LHS), for a large $K\in \N$, the decomposition: 
$$
F(n)=\sum_{d|n,d\le K}F'(d)+\sum_{d|n,d>K}F'(d) 
$$
\par
\noindent
rendering in the LHS the following (again, sums exchange is possible because $F'$ may not depend on $n$): 
$$
{1\over x}\sum_{n\le x}F(n)c_q(n)=\sum_{d\le K}F'(d){1\over x}\sum_{m\le x/d}c_q(dm)+\sum_{d>K}F'(d){1\over x}\sum_{m\le x/d}c_q(dm), 
$$
\par
\noindent
in which, now, we apply two different treatments, depending on $d\le K$ or $d>K$. For low divisors $d$, 
$$
\sum_{d\le K}F'(d){1\over x}\sum_{m\le x/d}c_q(dm)=\sum_{d\le K}F'(d)\sum_{j\le q,(j,q)=1}{1\over x}\sum_{m\le x/d}e_q(jdm)
$$
$$
=\sum_{d\le K}F'(d)\sum_{j\le q,(j,q)=1}\left({1\over d}\cdot \1_{d\equiv 0\bmod q}+O\left({1\over x}\left(1+{{\1_{d\not \equiv 0\bmod q}}\over {\left\Vert {{jd}\over q}\right\Vert}}\right)\right)\right)
=\varphi(q)\sum_{{d\le K}\atop {d\equiv 0\bmod q}}{{F'(d)}\over d}+O(1/x), 
$$
\par				
\noindent
from used-a-lot exponential sums cancellations, with a final $O-$constant not affecting the $x-$decay, while for high divisors $d$: 
$$
\sum_{d>K}F'(d){1\over x}\sum_{m\le x/d}c_q(dm)\ll \varphi(q)\sum_{d>K}{{|F'(d)|}\over d}, 
$$
\par
\noindent
uniformly in $x>0$, using the trivial bound $|c_q(n)|\le \varphi(q)$, $\forall n\in \Z$. In all, 
$$
{1\over x}\sum_{n\le x}F(n)c_q(n)=\varphi(q)\sum_{{d\le K}\atop {d\equiv 0\bmod q}}{{F'(d)}\over d}+O(1/x)+O\left(\varphi(q)\sum_{d>K}{{|F'(d)|}\over d}\right), 
$$
\par
\noindent
entailing 
$$
{1\over {\varphi(q)}}\lim_{x\to \infty}{1\over x}\sum_{n\le x}F(n)c_q(n)=\sum_{{d\le K}\atop {d\equiv 0\bmod q}}{{F'(d)}\over d}+O\left(\sum_{d>K}{{|F'(d)|}\over d}\right), 
$$
\par
\noindent
actually, giving the required equation, since from Delange Hypothesis the series $\sum_{d=1}^{\infty}{{|F'(d)|}\over d}$ converges, so errors in $O$ are infinitesimal with $K$, an arbitrarily large natural number (also, present LHS doesn't depend on it!). Last but not least, this also proves the convergence in RHS of these, say, $d\le K$-coeff.s (as $K\to \infty$). 
\par
\hfill QED (Wintner-Delange Formula) 
\smallskip
\par
Let's turn to the application of this Formula to our case $F(a)=C_{f,g_N}(N,a)$, getting that (since we are assuming $(DH)$ in hypotheses) we have the Carmichael formula, $(CF)$ above. Now (mimicking the proof of [CMu2] Theorem 1, $(ii)$ $\Rightarrow$ $(iii)$, exactly) we'll get the Reef above; in fact, let's calculate, since we know that the shift Ramanujan expansion converges (again, from $(DH)$ implying this by just proved Wintner-Delange), its shift-Ramanujan coefficients, for correlation $C_{f,g_N}(N,a)$, namely 
$$
\widehat{C_{f,g_N}}(N,\ell)={1\over {\varphi(\ell)}}\lim_{x\to \infty}{1\over x}\sum_{a\le x}C_{f,g_N}(N,a)c_{\ell}(a). 
$$
\par
\noindent
Plugging, so to speak, $(2)$ with $Q=N$ inside this RHS, we get for it : 
$$
{1\over x}\sum_{a\le x}C_{f,g_N}(N,a)c_{\ell}(a)=\sum_{q\le Q}\hat{g}(q)\sum_{n\le N}f(n){1\over x}\sum_{a\le x}c_q(n+a)c_{\ell}(a), 
$$
\par
\noindent
present exchange of sums being possible thanks to the hypothesis: $C_{f,g_N}(N,a)$ is fair. Then, 
$$
{1\over {\varphi(\ell)}}\lim_{x\to \infty}{1\over x}\sum_{a\le x}C_{f,g_N}(N,a)c_{\ell}(a)={1\over {\varphi(\ell)}}\sum_{q\le Q}\hat{g}(q)\sum_{n\le N}f(n)\lim_{x\to \infty}{1\over x}\sum_{a\le x}c_q(n+a)c_{\ell}(a), 
\leqno{(\ast \ast)}
$$
\par
\noindent
since all we are exchanging with ${\displaystyle \lim_{x\to \infty} }$ are finite sums (again, we're implicitly using fairness); then, the orthogonality of Ramanujan sums (first proved by Carmichael in [Ca], that's why $(CF)$ bears his name), namely Theorem 1 in [Mu]: 
$$
\lim_{x\to \infty}{1\over x}\sum_{a\le x}c_q(n+a)c_{\ell}(a)=\1_{q=\ell}\thinspace \cdot \thinspace c_q(n),
\enspace 
\forall \ell,n,q\in \N,
$$
\par
\noindent
gives inside $(\ast \ast)$ whence for quoted $(CF)$ the shift-Ramanujan coefficients 
$$
\widehat{C_{f,g_N}}(N,\ell)={1\over {\varphi(\ell)}}\hat{g}(\ell)\sum_{n\le N}f(n)c_{\ell}(n) 
$$
\par
\noindent
and this, thanks to the finite support of $\hat{g}$, up to $Q=N$, here, gives the R.e.e.f.! QED 
\smallskip
\par				
One last detail: equation $(2)$, actually, we didn't prove; but it follows from $m=n+a$ in (another unproven) 
$$
\sum_{q|m,q\le Q}g'(q)=\sum_{\ell \le Q}\hat{g}(\ell)c_{\ell}(m),
$$
\par
\noindent
that is : the $g_Q$ (see paper beginning) finite Ramanujan expansion, f.R.e. (for which we referred to [CMu1], of course), with Ramanujan coefficients
$$
\hat{g}(\ell)\defineq \sum_{q\equiv 0\bmod \ell}{{g'(q)}\over q}. 
$$
\par
\noindent
This can be proved at once, from quoted Lemma 1 of [CMu1], that we also prove (briefly) here: 
$$
\1_{q|m}={1\over q}\sum_{\ell|q}c_{\ell}(m), 
$$
\par
\noindent
because : the orthogonality of additive characters [Da] (rearranging by g.c.d.) gives 
$$
\1_{q|m}={1\over q}\sum_{r\le q}e_q(rm)
={1\over q}\sum_{\ell|q}\sum_{r\le q,(r,q)=q/\ell}e_q(rm)
={1\over q}\sum_{\ell|q}\sum_{j\le \ell,(j,\ell)=1}e_{\ell}(jm),
\enspace
\hbox{\rm with}
\enspace 
c_{\ell}(n)\defineq \sum_{{j\le \ell}\atop {(j,\ell)=1}}e_{\ell}(jn). 
$$
\par
\noindent
Then from this divisiblity condition we prove $g_Q$ f.R.e.: 
$$
\sum_{q|m,q\le Q}g'(q)=\sum_{q\le Q}{{g'(q)}\over q}\sum_{\ell|q}c_{\ell}(m)=\sum_{\ell \le Q}\hat{g}(\ell)c_{\ell}(m), 
$$
\par
\noindent
simply exchanging sums and using above definition of f.R.e. coefficients, $\hat{g}(q)$. QED (for equation $(2)$, too.)
\par
\hfill $\square$ 

\bigskip
\bigskip

\par
\centerline{\stampatello 3. The well-known case $f=g=\Lambda$, $a=2k>0$ of our Theorem : $2k-$prime-twins.} 
\smallskip
\par
\noindent
\par
\noindent
(Actually, in my talk I thought that the case we are exposing now could not be treated; but, taking $Q=N$ in Theorem 1 of [CMu2] and truncating $g$ as $g_N$ with the error in $(1)$, then, from this cut of original correlation $C_{f,g}=C_{\Lambda,\Lambda}$, the case of $2k-$twin primes is now contemplated. ) 
\smallskip
\par
Assuming $(DH)$ for $f=g=\Lambda$, Hardy-Littlewood heuristic (Conjecture B and (5.26) [HL]) is a Theorem. 

\medskip

\par
We apply, in fact, the calculations for Ramanujan coefficients of $N-$truncated von Mangoldt function, $\Lambda_N$, from the classical [Da] von Mangoldt $\Lambda=(-\mu \log)\ast \1$, [T], defined as usual in terms of primes $p\in \P$ : 
$$
\Lambda(n)\defineq \sum_{k\in \N}\sum_{p\in \P}\1_{n=p^k}\log p
\enspace 
\Rightarrow 
\enspace 
\Lambda(n)=\sum_{d|n}(-\mu(d)\log d),
\enspace 
\Lambda_N(n)=\sum_{d|n,d\le N}(-\mu(d)\log d), 
$$
\par
\noindent
entailing 
$$
\Lambda_N(n)=\sum_{q\le N}\widehat{\Lambda_N}(q)c_q(n),
\quad 
\widehat{\Lambda_N}(q)\defineq -\sum_{{d\le N}\atop {d\equiv 0\bmod q}}{{\mu(d)\log d}\over d}\ll {{\log^2 N}\over q}, 
$$
\par
\noindent
where now these are, thanks to $\S4$ of [CMu2], with an absolute $c>0$, 
$$
\widehat{\Lambda_N}(q)={{\mu(q)}\over {\varphi(q)}}+O\left({1\over q}\exp\left(-c\sqrt{\log N}\right)\right),
\enspace 
\forall q\le \sqrt{N}, 
$$
\par				
\noindent
thanks to the zero-free region of Riemann zeta-function (actually, we are not using most recent one). 
Now, 
$$
C_{\Lambda,\Lambda}(N,a)=\sum_{\ell \le N}{{\widehat{\Lambda_N}(\ell)}\over {\varphi(\ell)}}\left(\sum_{n\le N}\Lambda(n)c_{\ell}(n)\right)c_{\ell}(a)
 +O_{\varepsilon}\left(N^{\varepsilon}\left(N+a\right)^{\varepsilon}a\right), 
$$
\par
\noindent
from our Theorem: $C_{\Lambda,\Lambda_N}$ is fair \& assume $(DH)$, $f=g=\Lambda$; set \thinspace $a=2k>0$, ${{\log k}\over {\log N}}<1-\delta$, $\delta \in (0,1/2)$ fixed: 
$$
C_{\Lambda,\Lambda}(N,a)=\sum_{\ell \le \sqrt{N}}{{\mu(\ell)}\over {\varphi^2(\ell)}}\left(\sum_{n\le N}\Lambda(n)c_{\ell}(n)\right)c_{\ell}(a)
 + O\left(\exp\left(-c\sqrt{L}\right)\sum_{\ell \le \sqrt{N}}{{(a,\ell)}\over {\ell \varphi(\ell)}}\sum_{n\le N}\Lambda(n)(n,\ell)\right)
$$
$$
+O\left(L^2 \sum_{\sqrt{N}<\ell \le N}{{(a,\ell)}\over {\ell\varphi(\ell)}}\sum_{n\le N}\Lambda(n)(n,\ell)\right)+O\left(N^{1-\delta}\right), 
$$
\par
\noindent
where we have applied well-known $|c_q(n)|\le (q,n)$, see Lemma A.1 in [CMu2], and above bounds for $\Lambda_N$, abbreviating hereafter $L\defineq \log N$. In the main term, applying PNT(Prime Number Theorem) [Da], [T] : 
$$
\sum_{n\le N}\Lambda(n)c_{\ell}(n)=\mu(\ell)\sum_{{n\le N}\atop {(n,\ell)=1}}\Lambda(n)+O\Big(L\varphi(\ell)\sum_{p|\ell}\log p\Big)
\buildrel{{\rm PNT}}\over{=\!=}\mu(\ell)N+O\left(Ne^{-c\sqrt{L}}\right)+O\left(L\varphi(\ell)\log \ell\right), 
$$
\par
\noindent
from well known [Da]: \enspace $\sum_{p|\ell}\log p\le \sum_{n|\ell}\Lambda(n)=\log \ell$; here, we need to bound the $n-$sum in remainders as 
$$
\sum_{n\le N}\Lambda(n)(n,\ell)=\sum_{d|\ell}d\sum_{{n\le N}\atop {(n,\ell)=d}}\Lambda(n)
\ll \sum_{d|\ell}d\sum_{{n\le N}\atop {n\equiv 0\bmod d}}\Lambda(n)
\ll N+\ell L\sum_{k\in \N}\sum_{p^k|\ell}\log p
\ll NL^2, 
\enspace
\forall \ell \le N, 
$$
\par
\noindent
by \v{C}ebi\v{c}ev bound [T]: \enspace $\sum_{n\le N}\Lambda(n)\ll N$. Then, using [T]: \enspace $\varphi(\ell)\gg \ell/\log \ell$, changing time to time \enspace $c>0$, 
$$
C_{\Lambda,\Lambda}(N,a)=N\sum_{\ell \le \sqrt{N}}{{\mu^2(\ell)}\over {\varphi^2(\ell)}}c_{\ell}(a)
 +O\left(Ne^{-c\sqrt{L}}\sum_{\ell \le \sqrt{N}}{{(a,\ell)}\over {\ell^2}}+NL^5 \sum_{\sqrt{N}<\ell \le N}{{(a,\ell)}\over {\ell^2}}+N^{1-\delta}\right)
$$
$$
=N\sum_{\ell=1}^{\infty}{{\mu^2(\ell)}\over {\varphi^2(\ell)}}c_{\ell}(a)+O\left(N\sum_{\ell>\sqrt{N}}{{\log^2 \ell}\over {\ell^2}}(a,\ell)\right)
 +O\left(Ne^{-c\sqrt{L}}\sum_{\ell \le \sqrt{N}}{{(a,\ell)}\over {\ell^2}}+NL^5 \sum_{\sqrt{N}<\ell \le N}{{(a,\ell)}\over {\ell^2}}+N^{1-\delta}\right), 
$$
\par
\noindent
being, by the definition of classic {\it singular series} for $a=2k-$twin primes, \enspace 
$$
\SingSer_{\Lambda,\Lambda}(a)\defineq \sum_{\ell=1}^{\infty}{{\mu^2(\ell)}\over {\varphi^2(\ell)}}c_{\ell}(a)
$$
\par
\noindent
and, also, by following bounds: (use $(A+B)^2\ll A^2+B^2$, then, [T]: $\sum_{d|a}1\ll_{\varepsilon}a^{\varepsilon}$ and $\sum_{d\le x}1/d\ll \log x$) 
$$
\sum_{\ell>\sqrt{N}}{{\log^2 \ell}\over {\ell^2}}(a,\ell)\ll \sum_{{d|a}\atop {d\le \sqrt{N}}}{1\over d}\sum_{m>\sqrt{N}/d}{{\log^2 d+\log^2 m}\over {m^2}}+\sum_{{d|a}\atop {d>\sqrt{N}}}{1\over d}\sum_{m=1}^{\infty}{{\log^2 d+\log^2 m}\over {m^2}}
\ll_{\varepsilon}a^{\varepsilon}{{L^2}\over {\sqrt{N}}}, 
$$
$$
\sum_{\ell \le \sqrt{N}}{{(a,\ell)}\over {\ell^2}}\ll \sum_{{d|a}\atop {d\le \sqrt{N}}}{1\over d}\sum_{m\le \sqrt{N}/d}{1\over {m^2}}
\ll L,
$$
$$
\sum_{\sqrt{N}<\ell \le N}{{(a,\ell)}\over {\ell^2}}\ll 
 \sum_{{d|a}\atop {d\le \sqrt{N}}}{1\over d}\sum_{\sqrt{N}/d<m\le N/d}{1\over {m^2}}+\sum_{{d|a}\atop {d>\sqrt{N}}}{1\over d}\sum_{m\le N/d}{1\over {m^2}}
\ll_{\varepsilon}{{a^{\varepsilon}}\over {\sqrt{N}}},
$$
\par				
\noindent
uniformly in $a=2k$, $k\in \N$, with ${{\log k}\over {\log N}}<1-\delta$, for a fixed $\delta \in (0,1/2)$, proves Hardy-Littlewood Conjecture\footnote{$^4$}{In my talk's jargon, we reached the Reef, so this is our treasure !} 
$$
C_{\Lambda,\Lambda}(N,2k)=\SingSer_{\Lambda,\Lambda}(2k)N+O
\left
(
Ne^{-c\sqrt{\log N}}
\right
). 
$$
\medskip
\par
\noindent
We are sorry, we don't have time to deepen (but we've plenty of margins\footnote{$^5$}{In 1637 Fermat wrote \lq \lq ... Hanc marginis exiguitas non caperet.\rq \rq}). 

\medskip

I wish to thank Ram Murty, not only for the biggest part of the work laying behind present Theorem \& Corollary (coming, but not exclusively, from [CMu2] of course) but also for the real beginning, of my interest in Ramanujan expansions \& their applications to analytic number theory, thanks to his \lq \lq illuminating\rq \rq, say, survey [Mu] on \lq \lq Ramanujan series\rq \rq, ironically (in the good meaning)  leading to {\it finite} Ramanujan expansions! 

\bigskip

\par
\centerline{\stampatello References}
\medskip
\item{\bf [Ca]} Carmichael, R.D.\thinspace - \thinspace {\sl Expansions of arithmetical functions in infinite series} \thinspace - \thinspace Proc. London Math. Society {\bf 34} (1932), 1--26. $\underline{\tt MR\thinspace 1576142}$ 
\smallskip
\item{\bf [CMu1]} Coppola, G. and Murty, M.Ram and Saha, B.\thinspace - \thinspace {\sl Finite Ramanujan expansions and shifted convolution sums of arithmetical functions} \thinspace - \thinspace J. Number Theory {\bf 174} (2017), 78--92. 
\smallskip
\item{\bf [CMu2]} Coppola, G. and Murty, M.Ram\thinspace - \thinspace {\sl Finite Ramanujan expansions and shifted convolution sums of arithmetical functions, II} \thinspace - \thinspace {\tt arXiv}:1705.07193, to appear on JNT 
\smallskip
\item{\bf [Da]} \thinspace Davenport, H.\thinspace - \thinspace {\sl Multiplicative Number Theory} \thinspace - \thinspace Third Edition, GTM 74, Springer, New York, 2000. $\underline{\tt MR\enspace 2001f\!:\!11001}$
\smallskip
\item{\bf [De]} \thinspace Delange, H.\thinspace - \thinspace {\sl On Ramanujan expansions of certain arithmetical functions}, Acta Arith., {\bf 31}  (1976), 259--270. $\underline{\tt MR\enspace 432578}$ 
\smallskip
\item{\bf [HL]} Hardy, G.H. and Littlewood,  J.E. \thinspace - \thinspace {\sl Some problems of \lq \lq Partitio numerorum\rq \rq. III: On the expression of a number as a sum of primes} \thinspace - \thinspace Acta Math., {\bf 44} (1923), 1--70. 
\smallskip
\item{\bf [Mu]} \thinspace Murty, M.Ram\thinspace - \thinspace {\sl Ramanujan series for arithmetical functions}, Hardy-Ramanujan J., {\bf 36} (2013), 21-33. 
\smallskip
\item{\bf [ScSp]} Schwarz, W. and Spilker, J.\thinspace - \thinspace {\sl Arithmetical functions,  (An introduction to elementary and analytic properties of arithmetic functions and to some of their almost-periodic properties).} London Mathematical Society Lecture Note Series, {\bf 184}, Cambridge University Press, Cambridge, 1994. $\underline{\tt MR\enspace 1274248}$ 
\smallskip
\item{\bf [T]} \thinspace Tenenbaum, G.\thinspace - \thinspace {\sl Introduction to Analytic and Probabilistic Number Theory} \thinspace - \thinspace Cambridge Studies in Advanced Mathematics, {\bf 46}, Cambridge University Press, 1995. $\underline{\tt MR\enspace 97e\!:\!11005b}$
\smallskip
\item{\bf [W]} \thinspace Wintner, A.\thinspace - \thinspace {\sl Eratosthenian averages} \thinspace - \thinspace Waverly Press, Baltimore, MD, 1943. $\underline{\tt MR0015082}$ 

\bigskip
\bigskip
\bigskip

\leftline{\tt Giovanni Coppola}
\leftline{\tt Universit\`{a} degli Studi di Salerno}
\leftline{\tt Home address : Via Partenio 12 - 83100, Avellino (AV) - ITALY}
\leftline{\tt e-mail : giovanni.coppola@unina.it}
\leftline{\tt e-page : www.giovannicoppola.name}
\leftline{\tt e-site : www.researchgate.net}

\bye